# STATIONARY MAX-STABLE FIELDS ASSOCIATED TO NEGATIVE DEFINITE FUNCTIONS

By Zakhar Kabluchko, Martin Schlather and Laurens de Haan

*University of Göttingen, University of Göttingen and Erasmus University Rotterdam*

Let $W_i, i \in \mathbb{N}$, be independent copies of a zero-mean Gaussian process $\{W(t), t \in \mathbb{R}^d\}$ with stationary increments and variance $\sigma^2(t)$. Independently of $W_i$, let $\sum_{i=1}^\infty \delta_{U_i}$ be a Poisson point process on the real line with intensity $e^{-y} dy$. We show that the law of the random family of functions $\{V_i(\cdot), i \in \mathbb{N}\}$, where $V_i(t) = U_i + W_i(t) - \sigma^2(t)/2$, is translation invariant. In particular, the process $\eta(t) = \bigvee_{i=1}^\infty V_i(t)$ is a stationary max-stable process with standard Gumbel margins. The process $\eta$ arises as a limit of a suitably normalized and rescaled pointwise maximum of $n$ i.i.d. stationary Gaussian processes as $n \to \infty$ if and only if $W$ is a (nonisotropic) fractional Brownian motion on $\mathbb{R}^d$. Under suitable conditions on $W$, the process $\eta$ has a mixed moving maxima representation.

**1. Introduction.** A stochastic process $\{\eta(t), t \in \mathbb{R}^d\}$ is called *max-stable* if, for any $n \in \mathbb{N}$, the process $\{\bigvee_{k=1}^n \eta_k(t), t \in \mathbb{R}^d\}$ has the same law as $\{\eta(t) + \log n, t \in \mathbb{R}^d\}$, where $\eta_1, \ldots, \eta_n$ are independent copies of $\eta$. It follows from this definition that the marginal distributions of $\eta$ are of the form $\mathbb{P}[\eta(t) \leq x] = \exp(-e^{-x+b(t)})$ and, more generally, the finite-dimensional distributions of $\eta$ are multivariate max-stable distributions of Gumbel type [26]. Max-stable processes have been studied in [8, 10, 12, 16, 29] and [9], Part III. Note that it is common to consider max-stable processes with Fréchet (rather than Gumbel) marginals, so most authors work with the process $e^\eta$ instead of $\eta$.

A general description of *stationary* max-stable processes in terms of nonsingular flows on measure spaces was given in [12]. A usual approach to constructing examples of such processes is to use some sort of moving maxima (or, more generally, mixed moving maxima) representation; see [11, 14, 27, 31]. Another family of examples, based on stationary random processes, was given in [27]. Contrary to the mixed moving maxima processes, which were shown to be mixing in [30], this family entails a nonvanishing large-distance dependence within the max-stable process.

In this paper, we are mainly interested in a remarkable stationary max-stable process constructed by Brown and Resnick in [4]. Let us recall part of their result (see also Section 9.8 in [9] for the two-sided version given here).

THEOREM 1. *Let $W_i, i \in \mathbb{N}$, be independent copies of a standard Brownian motion $\{W(t), t \in \mathbb{R}\}$ and, independently of $W_i$, let $\sum_{i=1}^{\infty} \delta_{U_i}$ be a Poisson point process on $\mathbb{R}$ with intensity $e^{-y}\,dy$. Then, the process*

$$\eta(t) = \bigvee_{i=1}^{\infty} (U_i + W_i(t) - |t|/2) \tag{1}$$

*is a stationary max-stable process with standard Gumbel margins.*

A natural question arises as to whether further stationary max-stable processes can be constructed by replacing, in the above construction, the drifted Brownian motion $W(t) - |t|/2$ by other stochastic processes. Thus, we are interested in stochastic processes $\{\xi(t), t \in \mathbb{R}^d\}$ having the property that the process $\eta(t) = \bigvee_{i=1}^{\infty}(U_i + \xi_i(t))$ is stationary, where $U_i, i \in \mathbb{N}$, are as above and $\xi_i, i \in \mathbb{N}$, are independent copies of $\xi$. We call such processes $\xi$ *Brown–Resnick stationary*; see Section 2 for a more precise definition. In [4], two different proofs of Theorem 1 were given. One of them is based on the fact that $e^{-y}\,dy$ is an invariant measure for the Brownian motion with drift $-1/2$ and can be extended to show that some classes of processes with Markov property are Brown–Resnick stationary; see [5, 30]. The other proof, which uses the connection with the extreme-value theory of Gaussian processes, will be discussed later in Sections 6 and 8.

We are going to show that any Gaussian process with stationary increments becomes Brown–Resnick stationary after subtracting an appropriate drift term. Recall that a random process $\{W(t), t \in \mathbb{R}^d\}$ is said to have stationary increments if the law of $\{W(t+t_0) - W(t_0), t \in \mathbb{R}^d\}$ does not depend on the choice of $t_0 \in \mathbb{R}^d$. If $W$ is a Gaussian process with stationary increments (always supposed to have *zero mean*), then its law is completely characterized by what we shall call the *variogram*

$$\gamma(t) = \mathbb{E}(W(t+t_0) - W(t_0))^2, \qquad t \in \mathbb{R}^d,$$

and the variance $\sigma^2(t) = \operatorname{Var} W(t)$. It is well known that a function $\gamma: \mathbb{R}^d \to [0, \infty)$ with $\gamma(0) = 0$ is a variogram of some Gaussian process with stationary increments if and only if it is negative definite. The latter condition means that $\gamma(-t) = \gamma(t)$ for every $t \in \mathbb{R}^d$ and $\sum_{i,j=1}^{n} a_i a_j \gamma(t_i - t_j) \leq 0$ for every $t_1, \ldots, t_n \in \mathbb{R}^d$ and $a_1, \ldots, a_n \in \mathbb{R}$ satisfying $\sum_{i=1}^n a_i = 0$; see [2] for more on negative definite functions. Examples of Gaussian processes with stationary increments are provided by, for example, stationary Gaussian processes, their integrals (if $d=1$) and fractional (Lévy) Brownian motions, the latter being characterized by $W(0) = 0$ and $\gamma(t) = \|t\|^\alpha$ for some $\alpha \in (0, 2]$. Here, $\|t\|$ denotes the Euclidean norm of the vector $t$.

THEOREM 2. *Let $W_i, i \in \mathbb{N}$, be independent copies of a Gaussian process $\{W(t), t \in \mathbb{R}^d\}$ with stationary increments, variance $\sigma^2(t)$ and variogram $\gamma(t)$. Independently of $W_i$, let $\sum_{i=1}^\infty \delta_{U_i}$ be a Poisson point process on $\mathbb{R}$ with intensity $e^{-y}\,dy$. The process*

$$(2) \qquad \eta(t) = \bigvee_{i=1}^\infty (U_i + W_i(t) - \sigma^2(t)/2)$$

*is then a stationary max-stable process with standard Gumbel margins. The law of $\eta$ depends only on the variogram $\gamma$.*

The process $\eta$ defined above will be called the *Brown–Resnick process* associated to the variogram $\gamma$. If $\gamma(t) = |t|$ [i.e., if the corresponding Gaussian process $W$, under $W(0) = 0$, is a standard Brownian motion], then we recover the process of Theorem 1, originally considered in [4]. The Brown–Resnick process corresponding to the variogram $\gamma(t_1, \ldots, t_d) = |t_1| + \cdots + |t_d|$ was used as a model of extreme spatial rainfall in [6] and [13]. Another natural class of random processes, having the advantage of being isotropic, can be obtained by taking $\gamma(t) = \|t\|^\alpha$, $t \in \mathbb{R}^d$, for some $\alpha \in (0, 2]$. If $\alpha = 2$, the corresponding drifted Gaussian process $W(t) - \sigma^2(t)/2$ is a "random parabola" of the form $W(t) = \langle t, N \rangle - \|t\|^2/2$, where the random vector $N$ has the standard Gaussian distribution on $\mathbb{R}^d$ and we recover a process introduced in [15] and [17]; see also [11]. If $\gamma$ is bounded, then the process $W$ can be chosen to be stationary (after changing the variance and without changing the variogram; see, e.g., Proposition 7.13 in [2]) and $\eta$ belongs to the class of max-stable processes considered in Theorem 2 of [27].

Different Gaussian processes with stationary increments may have the same variogram. For example, let $\{W(t), t \in \mathbb{R}\}$ be a standard Brownian motion and let $f \in L^2(\mathbb{R})$. The process $W_f(t) = W(t) + \int_\mathbb{R} f(s)\,dW(s)$ then has the same variogram $\gamma(t) = |t|$ as $W$ and it is not difficult to see that the laws of $W_f$ and $W_g$ coincide if and only if $f = g$ a.s. The fact that different processes with the same variogram lead to the same $\eta$ is quite surprising, even in the particular case mentioned above.

The Brown–Resnick processes defined in Theorem 2 have no a priori connection to mixed moving maxima processes mentioned at the beginning of the paper. It was asked in [30] if the original Brown–Resnick process corresponding to $\gamma(t) = |t|$ has a representation as a mixed moving maxima process. We shall show in Section 5 that the answer is affirmative. More generally, it will be shown that the Brown–Resnick process corresponding to a Gaussian process $W$ with stationary increments has a mixed moving maxima representation provided that $\lim_{\|t\|\to\infty}(W(t) - \sigma^2(t)/2) = -\infty$ a.s.

The paper is organized as follows. In Section 2, we introduce the notion of Brown–Resnick stationarity. In Section 3, we prove a general criterion which allows one to decide whether a given random process $\xi$ has the property of Brown–Resnick stationarity in terms of the Laplace transform of the finite-dimensional distributions of $\xi$. This criterion is then used in Section 4 to prove Theorem 2. In Section 5, we show that Brown–Resnick processes of Theorem 2 have a mixed moving maxima representation under some conditions on the variogram $\gamma$. In Sections 6 and 7, we study, generalizing [4], extremes of a large number of independent Gaussian processes. An alternative proof of Theorem 2, in the case $W(0) = 0$, is given in Section 8.

REMARK 3. Two objects will appear frequently in our considerations: the Poisson point process $\sum_{i=1}^{\infty} \delta_{U_i}$ with intensity $e^{-y}\,dy$ on $\mathbb{R}$ and the standard *Gumbel distribution* $\exp(-e^{-y})$, which is the distribution of $\max_{i\in\mathbb{N}} U_i$. The transformation $y \mapsto e^y$ allows us to switch from Gumbel to the more common *Fréchet* notation. That is, if $Y$ is a random variable with standard Gumbel distribution, then $Z = e^Y$ has standard Fréchet distribution, meaning that $\mathbb{P}[Z \leq z] = \exp(-1/z)$, $z > 0$. Further, $\sum_{i=1}^{\infty} \delta_{e^{U_i}}$ is a Poisson point process on $(0, \infty)$ with intensity $dz/z^2$. Thus, if $\eta$ is a max-stable process, as defined at the beginning of the paper, then the process $e^\eta$ is max-stable in the usual sense [8].

**2. Brown–Resnick stationarity property.** Let $\xi_i, i \in \mathbb{N}$, be independent copies of a random process $\{\xi(t), t \in \mathbb{R}^d\}$ satisfying

(3)  $$\mathbb{E}e^{\xi(t)} < \infty \qquad \text{for all } t \in \mathbb{R}^d.$$

Further, let $\sum_{i=1}^{\infty} \delta_{U_i}$ be a Poisson point process on $\mathbb{R}$ with intensity $e^{-y}\,dy$, independent of the family $\xi_i$, $i \in \mathbb{N}$. Define a process $\{\eta(t), t \in \mathbb{R}^d\}$ by

(4)  $$\eta(t) = \bigvee_{i=1}^{\infty}(U_i + \xi_i(t)).$$

The process $\eta$ is necessarily max-stable [8]. To give a short proof of this fact, let $\eta_1, \ldots, \eta_n$ be independent copies of $\eta$, constructed by starting with $\sum_{i=1}^{\infty} \delta_{U_{i,k}}$, $k = 1, \ldots, n$, and $\xi_{i,k}$, $i \in \mathbb{N}$, $k = 1, \ldots, n$, all objects being independent. The superposition $\sum_{k=1}^{n} \sum_{i=1}^{\infty} \delta_{U_{i,k}}$ is then a Poisson point process on $\mathbb{R}$ with intensity $ne^{-y}\,dy = e^{-(y-\log n)}\,dy$. Hence, $\sum_{k=1}^{n} \sum_{i=1}^{\infty} \delta_{U_{i,k}-\log n}$ has the law of the Poisson point process with intensity $e^{-y}\,dy$. So, the process $\bigvee_{k=1}^{n} \eta_k - \log n$ has the same law as $\eta$, which proves the max-stability of $\eta$.

By [8], the converse is also true: any stochastically continuous max-stable process $\eta$ is of the form (4) for some process $\xi$. The finite-dimensional distributions of $\eta$ were computed in [8]: given $t_1, \ldots, t_n \in \mathbb{R}^d$ and $y_1, \ldots, y_n \in \mathbb{R}$, we have

$$\text{(5)} \qquad \mathbb{P}[\eta(t_1) \leq y_1, \ldots, \eta(t_n) \leq y_n] = \exp\Big\{-\mathbb{E}\exp\max_{i=1,\ldots,n}(\xi(t_i) - y_i)\Big\}.$$

In particular, condition (3) ensures that for every $t \in \mathbb{R}^d$, $\eta(t)$ is finite a.s. We are interested in processes $\xi$ leading to a *stationary* process $\eta$.

DEFINITION 4. A stochastic process $\{\xi(t), t \in \mathbb{R}^d\}$ satisfying (3) is called *Brown–Resnick stationary* if the process $\eta$ defined in (4) is stationary.

It is trivial that every stationary process satisfying (3) is Brown–Resnick stationary. However, the converse is not true: by a result of [4], the nonstationary process $\xi(t) = W(t) - |t|/2$, where $\{W(t), t \in \mathbb{R}\}$ is a standard Brownian motion, is Brown–Resnick stationary. The next proposition gives an equivalent, but perhaps more natural, version of Definition 4.

PROPOSITION 5. *A process $\{\xi(t), t \in \mathbb{R}^d\}$ which satisfies (3) is Brown–Resnick stationary if and only if $\sum_{i=1}^{\infty} \delta_{U_i + \xi_i(\cdot)}$ is a translation invariant Poisson point process on the space $E = \mathbb{R}^{\mathbb{R}^d}$.*

Before we can start the proof, we need to introduce some notation. We endow $E = \mathbb{R}^{\mathbb{R}^d}$, the space of real-valued functions on $\mathbb{R}^d$, with the product $\sigma$-algebra $\mathcal{B}(E)$ generated by the finite-dimensional cylinder sets, that is, by the sets of the form

$$\text{(6)} \qquad C_{t_1,\ldots,t_n}(B) = \{f : \mathbb{R}^d \to \mathbb{R} : (f(t_1), \ldots, f(t_n)) \in B\},$$

where $t_1, \ldots, t_n \in \mathbb{R}^d$ and $B$ is a Borel set in $\mathbb{R}^n$. If the processes $\xi_i$ have continuous sample paths, then $E = C(\mathbb{R}^d)$, the space of continuous functions, could be considered as well. Let $\mathcal{M}(E)$ be the space of all measures on $E$ which have the form $\mu = \sum_{i=1}^{\infty} \delta_{f_i}$ for some $f_i \in E$ and which are locally finite [i.e., finite on all cylinder sets of the form (6) with bounded $B$]. We endow $\mathcal{M}(E)$ with the $\sigma$-algebra $\mathcal{B}(\mathcal{M}(E))$ generated by the maps $F_{t_1,\ldots,t_n;B} : \mathcal{M}(E) \to \mathbb{N}_0 \cup \{\infty\}$, $\mu \mapsto \mu(C_{t_1,\ldots,t_n}(B))$. A point process on $E$ is a random variable $\Theta : \Omega \to \mathcal{M}(E)$, defined on some probability space $\Omega$ and taking values in $\mathcal{M}(E)$. Also, recall (see [19, 26]) that for a locally finite measure $\Lambda$ on $E$, a Poisson point process with intensity $\Lambda$ is a point process $\Theta : \Omega \to \mathcal{M}(E)$ such that $\Theta(\cdot)(A) \sim \text{Poiss}(\Lambda(A))$ for each $A \in \mathcal{B}(E)$, $\Lambda(A) < \infty$, and the random variables $\Theta(\cdot)(A_i)$, $i \in \mathbb{N}$, are independent provided $A_i \in \mathcal{B}(E)$ are disjoint.

We define a family of operators $T_h : \mathcal{M}(E) \to \mathcal{M}(E)$, $h \in \mathbb{R}^d$, as follows: for $\mu = \sum_{i=1}^\infty \delta_{f_i} \in \mathcal{M}(E)$, we define $T_h(\mu) = \sum_{i=1}^\infty \delta_{f_i(\cdot + h)}$. A point process on $E$ is called *translation invariant* if its distribution, viewed as a probability measure on $\mathcal{M}(E)$, is invariant with respect to the family $T_h$. A measure $\Lambda$ on the space $E$ is called translation invariant if, for every $A \in \mathcal{B}(E)$ and every $h \in \mathbb{R}^d$, we have $\Lambda(A) = \Lambda(\{f(\cdot + h) : f \in A\})$. A Poisson point process $\Theta$ on $E$ is translation invariant if and only if its intensity measure $\Lambda$ is translation invariant.

PROOF OF PROPOSITION 5.  Let $\mathbb{P}_\xi$ be the law of $\xi$ on the space $E = \mathbb{R}^{\mathbb{R}^d}$. Define a map $\pi : \mathbb{R} \times E \to E$ by $\pi(U, \xi(\cdot)) = U + \xi(\cdot)$ and let $\Lambda$ be the pushforward of the measure $e^{-y} dy \times d\mathbb{P}_\xi$ by the map $\pi$ [i.e., for $A \in \mathcal{B}(E)$, define $\Lambda(A) = \int_{\pi^{-1}(A)} e^{-y} dy \times d\mathbb{P}_\xi$]. We show that condition (3) implies that the measure $\Lambda$ is locally finite. To this end, take $t \in \mathbb{R}^d$ and let $A_{t,k} = \{f \in E : f(t) > k\}$, $k \in \mathbb{Z}$. Then,

$$\Lambda(A_{t,k}) = \int_\mathbb{R} e^{-y} \mathbb{P}[\xi(t) > k - y] dy = e^{-k} \int_\mathbb{R} e^z \mathbb{P}[\xi(t) > z] dz,$$

which is finite, by (3). Since any bounded cylinder set is contained in some $A_{t,k}$, the measure $\Lambda$ is locally finite. Since $\bigcup_{k \in \mathbb{Z}} A_k = E$, the measure $\Lambda$ is $\sigma$-finite.

The random measure $\sum_{i=1}^\infty \delta_{(U_i, \xi_i(\cdot))}$ may be viewed as a Poisson point process on $\mathbb{R} \times E$ with intensity $e^{-y} dy \times d\mathbb{P}_\xi$. Therefore, by a general mapping theorem (see [19]), $\sum_{i=1}^\infty \delta_{U_i + \xi_i(\cdot)}$ is a Poisson point process on $E$ with intensity measure $\Lambda$. Given $t_1, \ldots, t_n \in \mathbb{R}^d$, $y_1, \ldots, y_n \in \mathbb{R}$ and denoting $B = \mathbb{R}^n \setminus \bigtimes_{j=1}^n (-\infty, y_j]$, we have

$$\begin{aligned}
\mathbb{P}[\eta(t_1) \leq y_1, \ldots, \eta(t_n) \leq y_n] &= \mathbb{P}[\nexists i \in \mathbb{N} : U_i + \xi_i(\cdot) \in C_{t_1,\ldots,t_n}(B)] \\
&= \exp(-\Lambda(C_{t_1,\ldots,t_n}(B))).
\end{aligned} \quad (7)$$

Now, suppose that the point process $\sum_{i=1}^\infty \delta_{U_i + \xi_i(\cdot)}$ is translation invariant. It follows that its intensity measure $\Lambda$ is translation invariant. Equation (7) then implies that the process $\eta$ is stationary. Conversely, if $\eta$ is stationary, then, again using (7), we obtain that

$$\Lambda(C_{t_1+h,\ldots,t_n+h}(B)) = \Lambda(C_{t_1,\ldots,t_n}(B))$$

for every set $B$ of the form $\mathbb{R}^n \setminus \bigtimes_{j=1}^n (-\infty, y_j]$ and every $h \in \mathbb{R}^d$. The translation invariance of $\Lambda$ follows from this, using uniqueness of extension of measures and the $\sigma$-finiteness of $\Lambda$.  □

**3. A general stationarity criterion.** In this section, we prove a general criterion for the Brown–Resnick stationarity of a given process in terms of Laplace transforms of its finite-dimensional distributions. Let $\{\xi(t), t \in \mathbb{R}^d\}$ be a random process satisfying (3). For $t_1, \ldots, t_n \in \mathbb{R}^d$, denote by $\mathbb{P}_{t_1,\ldots,t_n}$ the distribution of the random vector $(\xi(t_1), \ldots, \xi(t_n))$. An application of Hölder's inequality shows that the Laplace transform of the measure $\mathbb{P}_{t_1,\ldots,t_n}$, defined by

$$\varphi_{t_1,\ldots,t_n}(u_1, \ldots, u_n) = \int_{\mathbb{R}^n} e^{u_1 x_1 + \cdots + u_n x_n} \, d\mathbb{P}_{t_1,\ldots,t_n}(x_1, \ldots, x_n),$$

is finite provided $u_i \in [0,1]$, $\sum_{i=1}^n u_i \leq 1$.

PROPOSITION 6. *A random process $\{\xi(t), t \in \mathbb{R}^d\}$ satisfying the moment condition (3) is Brown–Resnick stationary if and only if*

(8) $$\varphi_{t_1,\ldots,t_n}(u_1, \ldots, u_n) = \varphi_{t_1+h,\ldots,t_n+h}(u_1, \ldots, u_n)$$

*for every $h, t_1, \ldots, t_n \in \mathbb{R}^d$ and any $u_1, \ldots, u_n \in [0,1]$ satisfying $\sum_{i=1}^n u_i = 1$.*

We need the following lemma on the uniqueness of the Laplace transform.

LEMMA 7. *Let $\mu_1$ and $\mu_2$ be two finite measures on $\mathbb{R}^n$ with Laplace transforms $\psi_1(t) = \int_{\mathbb{R}^n} e^{\langle t,s\rangle} \, d\mu_1(s)$ and $\psi_2(t) = \int_{\mathbb{R}^n} e^{\langle t,s\rangle} \, d\mu_2(s)$ such that $\psi_1$ and $\psi_2$ are finite and equal on some open set $D \subset \mathbb{R}^n$. Then, $\mu_1 = \mu_2$.*

PROOF. If $\psi_1$ and $\psi_2$ are finite on $D$, then they are finite on the complexification of $D$, that is, on the set $D^c = \{t \in \mathbb{C}^n : \operatorname{Re} t \in D\}$. Since $\psi_1$ and $\psi_2$ are analytic functions coinciding on $D$, they must coincide on $D^c$. Let $t_0 \in D$. Then, $s \mapsto \psi_1(t_0 + is)$ is the characteristic function of the finite measure $e^{\langle t_0, \cdot\rangle} \, d\mu_1(\cdot)$. Now, $\psi_1(t_0 + is) = \psi_2(t_0 + is)$ and the fact that a finite measure is uniquely determined by its characteristic function together imply that $e^{\langle t_0, \cdot\rangle} \, d\mu_1(\cdot) = e^{\langle t_0, \cdot\rangle} \, d\mu_2(\cdot)$. Hence, $\mu_1 = \mu_2$. □

PROOF OF PROPOSITION 6. We use the notation of the previous section. Our goal is to show that the intensity measure $\Lambda$ is translation invariant if and only if (8) holds. For a set $B \subset \mathbb{R}^n$ and $x \in \mathbb{R}$, let $B + x = B + (x, x, \ldots, x)$. For a cylinder set $C_{t_1,\ldots,t_n}(B)$ [recall (6)], we have

$$\Lambda(C_{t_1,\ldots,t_n}(B))$$
$$= \int_{\mathbb{R}} e^x \mathbb{P}_{t_1,\ldots,t_n}(B + x) \, dx$$
$$= \int_{\mathbb{R}} \int_{\mathbb{R}^n} e^x 1_{B+x}(y_1, \ldots, y_n) \, d\mathbb{P}_{t_1,\ldots,t_n}(y_1, \ldots, y_n) \, dx$$

$$= \int_{\mathbb{R}} \int_{\mathbb{R}^n} e^{y_1} e^{x-y_1} 1_{B+x-y_1}(0, y_2 - y_1, \ldots, y_n - y_1) \, d\mathbb{P}_{t_1,\ldots,t_n}(y_1, \ldots, y_n) \, dx$$

$$= \int_{\mathbb{R}} \int_{\mathbb{R}^n} e^{y_1} e^{z} 1_{B+z}(0, y_2 - y_1, \ldots, y_n - y_1) \, d\mathbb{P}_{t_1,\ldots,t_n}(y_1, \ldots, y_n) \, dz.$$

Consider a measure $\mu_{t_1,\ldots,t_n}$ on $\mathbb{R}^n$, defined on Borel sets $A \subset \mathbb{R}^n$ by

$$\mu_{t_1,\ldots,t_n}(A) = \int_{\mathbb{R}^n} e^{y_1} 1_A(0, y_2 - y_1, \ldots, y_n - y_1) \, d\mathbb{P}_{t_1,\ldots,t_n}(y_1, \ldots, y_n).$$

Note that, by (3), we have $\mu_{t_1,\ldots,t_n}(A) \leq \mathbb{E}e^{\xi(t_1)} < \infty$ and therefore the measure $\mu_{t_1,\ldots,t_n}$ is finite. The measure $\mu_{t_1,\ldots,t_n}$ may be viewed as a type of exponentially weighted projection of the measure $\mathbb{P}_{t_1,\ldots,t_n}$ onto the $(n-1)$-dimensional hyperplane $\{(x_i)_{i=1}^n \in \mathbb{R}^n : x_1 = 0\}$. We have

(9) $$\Lambda(C_{t_1,\ldots,t_n}(B)) = \int_{\mathbb{R}} e^z \mu_{t_1,\ldots,t_n}(B + z) \, dz.$$

The Laplace transform of $\mu_{t_1,\ldots,t_n}$ is given by

(10) $$\begin{aligned} \psi_{t_1,\ldots,t_n}(u_1, \ldots, u_n) \\ &= \int_{\mathbb{R}^n} e^{y_1} e^{u_2(y_2-y_1)+\cdots+u_n(y_n-y_1)} \, d\mathbb{P}_{t_1,\ldots,t_n}(y_1, \ldots, y_n) \\ &= \int_{\mathbb{R}^n} e^{y_1(1-\sum_{i=2}^n u_i)+y_2 u_2+\cdots+y_n u_n} \, d\mathbb{P}_{t_1,\ldots,t_n}(y_1, \ldots, y_n) \\ &= \varphi_{t_1,\ldots,t_n}\left(1 - \sum_{i=2}^n u_i, u_2, \ldots, u_n\right), \end{aligned}$$

where $\varphi_{t_1,\ldots,t_n}$ is the Laplace transform of the measure $\mathbb{P}_{t_1,\ldots,t_n}$. Note that $\psi_{t_1,\ldots,t_n}$ does not depend on $u_1$.

Now, suppose that (8) holds. We then obtain

(11) $$\psi_{t_1,\ldots,t_n}(u_1, \ldots, u_n) = \psi_{t_1+h,\ldots,t_n+h}(u_1, \ldots, u_n)$$

provided that $u_i \in [0,1]$, $\sum_{i=2}^n u_i \leq 1$, which, by Lemma 7, implies that $\mu_{t_1,\ldots,t_n} = \mu_{t_1+h,\ldots,t_n+h}$ and hence, by (9),

(12) $$\Lambda(C_{t_1+h,\ldots,t_n+h}(B)) = \Lambda(C_{t_1,\ldots,t_n}(B)).$$

This proves the translation invariance of $\Lambda$ on the semi-ring of the cylinder sets. Using the theorem on the uniqueness of the extension of measures and the fact that $\Lambda$ is $\sigma$-finite, we obtain the translation invariance of $\Lambda$ on the whole $\sigma$-algebra $\mathcal{B}(E)$.

Now, suppose that $\Lambda$ is translation invariant. It follows that (12) holds and thus, using (9),

$$\int_{\mathbb{R}} e^z \mu_{t_1,\ldots,t_n}(B+z) \, dz = \int_{\mathbb{R}} e^z \mu_{t_1+h,\ldots,t_n+h}(B+z) \, dz$$

for every Borel set $B \subset \mathbb{R}^n$ and every $h, t_1, \ldots, t_n \in \mathbb{R}^d$. Since the measure $\mu_{t_1,\ldots,t_n}$ is concentrated on the hyperplane $\{(x_i)_{i=1}^n \in \mathbb{R}^n : x_1 = 0\}$, it follows that, actually, $\mu_{t_1,\ldots,t_n} = \mu_{t_1+h,\ldots,t_n+h}$. By considering the Laplace transforms, we obtain that (11) holds, from which (8) follows. This completes the proof. □

As an immediate consequence of the above proposition, we obtain the following nontrivial corollaries:

COROLLARY 8. *Let $\{\xi'(t), t \in \mathbb{R}^d\}$ and $\{\xi''(t), t \in \mathbb{R}^d\}$ be two independent processes, both having the Brown–Resnick stationarity property. The process $\xi' + \xi''$ is then also Brown–Resnick stationary.*

COROLLARY 9. *Let $\{\xi_1(t), t \in \mathbb{R}^{d_1}\}$ and $\{\xi_2(t), t \in \mathbb{R}^{d_2}\}$ be independent Brown–Resnick stationary processes. The process $\{\xi(t_1, t_2), (t_1, t_2) \in \mathbb{R}^{d_1+d_2}\}$ defined by $\xi(t_1, t_2) = \xi_1(t_1) + \xi_2(t_2)$ is then Brown–Resnick stationary.*

## 4. Max-stable processes associated to variograms.

THEOREM 10. *Let $\{W(t), t \in \mathbb{R}^d\}$ be a Gaussian process with stationary increments and variance $\sigma^2(t)$. The process $\xi(t) = W(t) - \sigma^2(t)/2$ is then Brown–Resnick stationary.*

PROOF. Recall our standing assumption $\mathbb{E}(W(t)) = 0$. It follows from the definition of the variogram $\gamma(t) = \mathbb{E}(W(t) - W(0))^2$ that we have

$$\operatorname{Cov}(W(t), W(s)) = \sigma^2(t)/2 + \sigma^2(s)/2 - \gamma(t-s)/2.$$

We are going to apply Proposition 6 to $\xi(t)$. Note that $\mathbb{E}e^{\xi(t)} = 1$, which shows that (3) is satisfied. We need to prove that (8) holds. The distribution $\mathbb{P}_{t_1,\ldots,t_n}$ of the random vector $(\xi(t_1), \ldots, \xi(t_n))$ is a multivariate Gaussian distribution whose expectation vector $(\mu_i)_{i=1,\ldots,n}$ and covariance matrix $(\sigma_{ij})_{i,j=1,\ldots,n}$ are given, respectively, by

(13) $\quad \mu_i = -\sigma^2(t_i)/2 \quad \text{and} \quad \sigma_{ij} = \sigma^2(t_i)/2 + \sigma^2(t_j)/2 - \gamma(t_i - t_j)/2.$

The Laplace transform of $\mathbb{P}_{t_1,\ldots,t_n}$ is given by

(14) $\quad \varphi_{t_1,\ldots,t_n}(u_1, \ldots, u_n) = \exp\left(\sum_{i=1}^n \mu_i u_i + \frac{1}{2} \sum_{i,j=1}^n \sigma_{ij} u_i u_j\right).$

Let $u_1, \ldots, u_n \in [0,1]$ satisfy $\sum_{i=1}^n u_i = 1$. By substituting $u_1 = 1 - \sum_{i=2}^n u_i$ into (14) and using (13), we obtain that

(15) $\quad \varphi_{t_1,\ldots,t_n}(u_1, \ldots, u_n) = \exp(L + \tfrac{1}{2}Q),$

where $L = L_{t_1,\ldots,t_n}(u_2,\ldots,u_n)$ and $Q = Q_{t_1,\ldots,t_n}(u_2,\ldots,u_n)$ are the linear part and the quadratic part, respectively (the constant term is easily seen to be zero). The linear part is given by

$$(16) \qquad L = \sum_{i=2}^{n}(\mu_i - \mu_1 + \sigma_{1i} - \sigma_{11})u_i = -\frac{1}{2}\sum_{i=2}^{n}\gamma(t_i - t_1)u_i.$$

The quadratic part is easily seen to be

$$(17) \qquad \begin{aligned} Q &= \sum_{i,j=2}^{n}(\sigma_{ij} - \sigma_{1i} - \sigma_{1j} + \sigma_{11})u_i u_j \\ &= \frac{1}{2}\sum_{i,j=2}^{n}(\gamma(t_i - t_1) + \gamma(t_j - t_1) - \gamma(t_j - t_i))u_i u_j. \end{aligned}$$

Thus, both terms $L$ and $Q$ do not change if one replaces $t_1,\ldots,t_n$ by $t_1 + h,\ldots,t_n + h$. Hence, (8) holds and the proof is complete. $\square$

PROPOSITION 11. *Let $W'$ and $W''$ be two Gaussian processes with stationary increments, having the same variogram $\gamma(t)$ and possibly different variances $\sigma'^2(t)$ and $\sigma''^2(t)$. Let $\Lambda'$ (resp., $\Lambda''$) be the intensity of the Poisson point process constructed as in Section 2 with $\xi$ replaced by $W' - \sigma'^2/2$ (resp., $W'' - \sigma''^2/2$). Then, $\Lambda' = \Lambda''$.*

PROOF. Formulas (15), (16) and (17) of the previous proof show that $\varphi'_{t_1,\ldots,t_n} = \varphi''_{t_1,\ldots,t_n}$, which, by (10), implies that $\psi'_{t_1,\ldots,t_n} = \psi''_{t_1,\ldots,t_n}$. Here, all objects marked with $'$ (resp., $''$) correspond to $W'$ (resp., $W''$). Lemma 7 yields $\mu'_{t_1,\ldots,t_n} = \mu''_{t_1,\ldots,t_n}$. Now, (9) shows that for every cylinder set $C_{t_1,\ldots,t_n}(B)$, we have

$$\Lambda'(C_{t_1,\ldots,t_n}(B)) = \Lambda''(C_{t_1,\ldots,t_n}(B)).$$

To finish the proof, use the $\sigma$-finiteness of $\Lambda'$ and $\Lambda''$. $\square$

REMARK 12. Given a Gaussian process $W$ with stationary increments, it will often be convenient to replace it by the process $\tilde{W}(t) = W(t) - W(0)$ having the same variogram $\gamma$ as $W$ and $\tilde{W}(0) = 0$. Note that the variance of the process $\tilde{W}$ is given by $\tilde{\sigma}^2(t) = \gamma(t)$.

PROOF OF THEOREM 2. The stationarity of $\eta$ follows from Theorem 10, whereas the max-stability was proven in the discussion following (4). The fact that $\eta(t)$ is standard Gumbel for each $t \in \mathbb{R}^d$ follows from (5). Finally, the last claim of the theorem follows from Proposition 11. $\square$

PROPOSITION 13. *If all Gaussian processes in Theorem 2 have continuous sample paths, then the process $\eta$ is also sample-continuous.*

PROOF. Let $K \subset \mathbb{R}^d$ be bounded. We use the notation $\xi(t) = W(t) - \sigma^2(t)/2$ and $\xi_i(t) = W_i(t) - \sigma^2(t)/2$. First, we show that for every $k \in \mathbb{Z}$, the random set
$$I_k = \Big\{ i \in \mathbb{N} : \sup_{t \in K}(U_i + \xi_i(t)) > k \Big\}$$
is a.s. finite. Indeed, the cardinality of $I_k$ is Poisson distributed with some (maybe infinite) intensity $\lambda_k$. We have
$$\lambda_k = \int_{-\infty}^{\infty} e^{-z} \mathbb{P}\Big[\sup_{t \in K} \xi(t) > k - z\Big] dz \leq 1 + \int_0^{\infty} e^z \mathbb{P}\Big[\sup_{t \in K} \xi(t) > k + z\Big] dz.$$
Since the process $\xi$ is Gaussian with continuous paths, a result of [20] (or see Corollary 3.2 of [22]) states that $\mathbb{E} \exp\{\varepsilon(\sup_{t \in K} \xi(t))^2\} < \infty$ for some small $\varepsilon > 0$. Hence, $\lambda_k < \infty$ and, consequently, $I_k$ is finite a.s.

We now show that $\eta$ is continuous a.s. Let $A_k$, $k \in \mathbb{Z}$, be the random event $\inf_{t \in K}(U_1 + \xi_1(t)) > k$. Note that $\mathbb{P}[\bigcup_{k \in \mathbb{Z}} A_k] = 1$. If, say, $A_k$ occurs, then
$$\eta(t) = \bigvee_{i \in I_k \cup \{1\}} (U_i + \xi_i(t)), \qquad t \in K.$$
It follows that $\eta$, being a pointwise maximum of a *finite* number of continuous functions, is itself continuous. □

**5. Representation as mixed moving maxima process.** We are now going to show that under some conditions on the underlying variogram $\gamma$, the Brown–Resnick process $\eta$ has a representation as a mixed moving maxima process. First, we recall a definition of mixed moving maxima processes as given in [27, 30]; see also [14, 29, 31]. Let $\{F(t), t \in \mathbb{R}^d\}$ be a measurable process and suppose that $\mathbb{E} \int_{\mathbb{R}^d} e^{F(t)} dt < \infty$. Let $\sum_{i=1}^{\infty} \delta_{(t_i, y_i)}$ be a Poisson point process on $\mathbb{R}^d \times \mathbb{R}$ with intensity $e^{-y} dt\, dy$ ($dt$ is the Lebesgue measure on $\mathbb{R}^d$) and let $F_i, i \in \mathbb{N}$, be independent copies of $F$. A process of the form
$$\eta(t) = \bigvee_{i=1}^{\infty} (F_i(t - t_i) + y_i), \qquad t \in \mathbb{R}^d,$$
is called a *mixed moving maxima process*. It is convenient to think of $F_i$ as a random mark attached to the point $(t_i, y_i)$. The process $\eta$ is stationary and max-stable; its finite-dimensional distributions are given by
$$\mathbb{P}[\eta(s_1) \leq z_1, \ldots, \eta(s_n) \leq z_n] = \exp\Big\{-\mathbb{E} \int_{\mathbb{R}^d} \exp \max_{j=1,\ldots,n}(F(s_j - t) - z_j)\, dt\Big\},$$
where $s_1, \ldots, s_n \in \mathbb{R}^d$, $z_1, \ldots, z_n \in \mathbb{R}$ and $\mathbb{E}$ denotes the expectation with respect to the law of $F$ (see, e.g., [29]).

THEOREM 14. *Let $\{W(t), t \in \mathbb{R}^d\}$ be a sample-continuous Gaussian process with stationary increments and variance $\sigma^2(t)$. Suppose that*

$$(18) \qquad \lim_{\|t\| \to \infty} (W(t) - \sigma^2(t)/2) = -\infty \qquad a.s.$$

*The Brown–Resnick process $\eta$ defined in Theorem 2 then has a representation as a mixed moving maxima process.*

PROOF. Recall that $\sum_{i=1}^{\infty} \delta_{U_i}$ is a Poisson point process on $\mathbb{R}$ with intensity $e^{-y} dy$ and $W_i, i \in \mathbb{N}$, are independent copies of $W$. The idea of the proof is to look at the random path $W_i(t) - \sigma^2(t)/2$, not from its *starting point* corresponding to $t = 0$, but rather from its *top point*. Let us be more precise.

Condition (18) implies that we may define a triple $(T, M, F) \in \mathbb{R}^d \times \mathbb{R} \times C(\mathbb{R}^d)$ by $M = \sup_{t \in \mathbb{R}^d}(W(t) - \sigma^2(t)/2)$, $T = \inf\{t \in \mathbb{R}^d : W(t) - \sigma^2(t)/2 = M\}$ (the "inf" is understood in, e.g., the lexicographic sense) and $F(t) = W(t+T) - \sigma^2(t+T)/2 - M$. So, $(T, M)$ are the coordinates of the top of the path $W(t) - \sigma^2(t)/2$, whereas $F(t)$ is the path itself, as viewed from its top. Let $M_i$, $T_i$ and $F_i$ be defined analogously, with $W$ replaced by $W_i$. Define a measurable transformation

$$\pi : \mathbb{R} \times C(\mathbb{R}^d) \to \mathbb{R}^d \times \mathbb{R} \times C(\mathbb{R}^d)$$

by mapping a pair $(U, W) \in \mathbb{R} \times C(\mathbb{R}^d)$ to the triple $(T, U + M, F) \in \mathbb{R}^d \times \mathbb{R} \times C(\mathbb{R}^d)$. Note that $\sum_{i=1}^{\infty} \delta_{(U_i, W_i)}$ is a Poisson point process on $\mathbb{R} \times C(\mathbb{R}^d)$ with intensity $e^{-y} dy \times d\mathbb{P}_W$, where $\mathbb{P}_W$ is the law of $W$ on $C(\mathbb{R}^d)$. Therefore, by the mapping theorem for Poisson point processes (see, e.g., [19]), we obtain that $\sum_{i=1}^{\infty} \delta_{(T_i, U_i + M_i, F_i)}$ is a Poisson point process on $\mathbb{R}^d \times \mathbb{R} \times C(\mathbb{R}^d)$ whose intensity measure $\Psi$ is given by

$$(19) \quad \Psi(A) = \int_{\pi^{-1}(A)} e^{-y} dy \times d\mathbb{P}_W = \int_{\mathbb{R}} e^{-y} \mathbb{P}[(T, M + y, F) \in A] dy,$$

where $A$ denotes a Borel subset of $\mathbb{R}^d \times \mathbb{R} \times C(\mathbb{R}^d)$.

We claim that the measure $\Psi$ has natural invariance properties. First, it follows from (19) that for every $z \in \mathbb{R}$, we have

$$\begin{aligned}
\Psi(A + (0, z, 0)) &= \int_{\mathbb{R}} e^{-y} \mathbb{P}[(T, M + y, F) \in A + (0, z, 0)] dy \\
&= \int_{\mathbb{R}} e^{-y} \mathbb{P}[(T, M + (y - z), F) \in A] dy \\
&= \int_{\mathbb{R}} e^{-(y+z)} \mathbb{P}[(T, M + y, F) \in A] dy \\
&= e^{-z} \Psi(A).
\end{aligned}$$

Second, Theorem 10 and Proposition 5 imply that $\Psi(A+(t,0,0)) = \Psi(A)$ for every $t \in \mathbb{R}^d$. To see this, note that the collection $\{(T_i, U_i + M_i, F_i), i \in \mathbb{N}\}$ can be obtained from the collection $\{U_i + W_i(\cdot) - \sigma^2(\cdot), i \in \mathbb{N}\}$, viewed as a translation invariant Poisson point process on $C(\mathbb{R}^d)$, by a measurable transformation, which commutes with spatial translations. Furthermore, note that

$$\Psi([0,1]^d \times [0,1] \times C(\mathbb{R}^d)) \leq \int_\mathbb{R} e^{-y} \mathbb{P}\Big[\sup_{t \in [0,1]^d} W(t) \geq -y\Big] dy$$

is finite by the same argument (based on [20]) as in the proof of Proposition 13.

We now show that the above invariance properties imply a product-type representation for $\Psi$. Take a measurable set $A \subset C(\mathbb{R}^d)$ and consider a measure $\Psi_A$ on $\mathbb{R}^d \times \mathbb{R}$, defined as follows: for $B \subset \mathbb{R}^d \times \mathbb{R}$, we set $\Psi_A(B) = \int_{B \times A} e^y d\Psi(t, y, F)$. By the above, the measure $\Psi_A$ is translation invariant and $\Psi_A([0,1]^d \times [0,1]) < \infty$. It follows that $\Psi_A$ is a multiple of the Lebesgue measure and hence we may write $d\Psi_A = \mathbb{Q}(A) \, dt \, dy$ for some finite constant $\mathbb{Q}(A)$. Further, $A \mapsto \mathbb{Q}(A)$ defines a finite measure on $C(\mathbb{R}^d)$. Introducing the normalized measure $\mathbb{Q}' = \mathbb{Q}/c$, where $c = \mathbb{Q}(C(\mathbb{R}^d))$, we may write $d\Psi$ in the form $ce^{-y} dt\, dy \times d\mathbb{Q}'$.

We are ready to finish the proof. The Brown–Resnick process of Theorem 2 may be written as

$$\eta(t) = \bigvee_{i=1}^\infty (U_i + W_i(t) - \sigma^2(t)/2) = \bigvee_{i=1}^\infty (F_i^*(t - t_i^*) + y_i^*),$$

where $F_i^*(\cdot) = F_i(\cdot) + \log c$, $t_i^* = T_i$ and $y_i^* = U_i + M_i - \log c$. We claim that this gives the required mixed moving maxima representation of $\eta$. First, recall that $\sum_{i=1}^\infty \delta_{(T_i, U_i+M_i, F_i)}$ is a Poisson point process on $\mathbb{R}^d \times \mathbb{R} \times C(\mathbb{R}^d)$ with intensity $d\Psi = ce^{-y} dt\, dy \times d\mathbb{Q}'$. It follows that $\sum_{i=1}^\infty \delta_{(t_i^*, y_i^*, F_i^*)}$ is a Poisson point process on the same space with intensity $e^{-y} dt\, dy \times d\mathbb{Q}^*$, where $\mathbb{Q}^*$ is the law of $F + \log c$ for $F \sim \mathbb{Q}'$. Thus, $\sum_{i=1}^\infty \delta_{(t_i^*, y_i^*)}$ is a Poisson point process on $\mathbb{R}^d \times \mathbb{R}$ with intensity $e^{-y} dt\, dy$, whereas $F_i^*$ may be viewed as a random mark sampled independently of $(t_i^*, y_i^*)$ according to the probability measure $\mathbb{Q}^*$, as required. □

REMARK 15. In the case $d=1$, it follows from Corollary 2.4 of [23] that condition (18) is satisfied whenever $\liminf_{t \to \infty} \gamma(t)/\log t > 8$.

**6. Maxima of independent Gaussian processes.** It was shown by Brown and Resnick [4] that a suitably normalized and spatially rescaled maximum of $n$ independent Brownian motions or Ornstein–Uhlenbeck processes converges, as $n \to \infty$, to the process $\eta$ of Theorem 1. Some related results were obtained in [15, 17, 18, 24]. We are going to extend the result of [4] to Gaussian processes whose covariance function satisfies a natural regular variation condition.

ASSUMPTION 16. Let $\{X(t), t \in D\}$ be a zero-mean, unit-variance Gaussian process defined on a neighborhood $D \subset \mathbb{R}^d$ of 0 and having covariance function $C(t_1, t_2) = \mathbb{E}[X(t_1)X(t_2)]$. We assume that the asymptotic relation

$$\lim_{\varepsilon \searrow 0} \frac{1 - C(\varepsilon t_1, \varepsilon t_2)}{L(\varepsilon)\varepsilon^\alpha} = \gamma(t_1 - t_2) \tag{20}$$

holds uniformly in $t_1, t_2 \in \mathbb{R}^d$ as long as $t_1, t_2$ stay bounded. Here, $L$ is a function varying slowly at 0, $\alpha \in (0, 2]$, and $\gamma : \mathbb{R}^d \to [0, \infty)$ is a continuous function satisfying $\gamma(\lambda t) = \lambda^\alpha \gamma(t)$ for every $\lambda \geq 0$, $t \in \mathbb{R}^d$.

Define normalizing sequences

$$b_n = (2 \log n)^{1/2} - (2 \log n)^{-1/2}((1/2) \log \log n + \log(2\sqrt{\pi})), \tag{21}$$

$$s_n = \inf\{s > 0 : L(s)s^\alpha = b_n^{-2}\} \tag{22}$$

and recall (see, e.g., Theorem 1.5.3 in [21]) that, for i.i.d. standard Gaussian $\{Z_i, i \in \mathbb{N}\}$, we have

$$\lim_{n \to \infty} \mathbb{P}\left[\bigvee_{i=1}^n b_n(Z_i - b_n) \leq y\right] = \exp(-e^{-y}). \tag{23}$$

We write $\eta_n \Rightarrow \eta$ as $n \to \infty$ if, for every compact set $K \subset \mathbb{R}^d$, the sequence of stochastic processes $\eta_n$ converges to $\eta$ weakly on $C(K)$, the space of continuous functions on $K$.

THEOREM 17. Let $X_i, i \in \mathbb{N}$, be independent sample-continuous copies of $X$, a process satisfying Assumption 16. Define

$$\eta_n(t) = \bigvee_{i=1}^n b_n(X_i(s_n t) - b_n).$$

Then, $\eta_n \Rightarrow \eta$ as $n \to \infty$, where $\eta$ is the Brown–Resnick process associated to the variogram $2\gamma$. In particular, $\gamma$ must be a variogram.

REMARK 18. The results of [4] can be recovered by applying the above theorem to the Ornstein–Uhlenbeck process and to the process $X(t) = B(t_0 + t)/(t_0 + t)^{1/2}$, where $t_0 > 0$ and $\{B(t), t \in \mathbb{R}\}$ is a standard Brownian motion.

PROOF OF THEOREM 17. Note that $s_n \to 0$ as $n \to \infty$. Define a process
$$Y_n(t) = b_n(X(s_n t) - b_n), \qquad t \in s_n^{-1} D.$$

Further, for $w \in \mathbb{R}$, let $Y_n^w$ be the process $Y_n - w$ conditioned on $\{Y_n(0) = w\}$. Let $Y_{i,n}$ and $Y_{i,n}^w$ be defined analogously, with $X$ replaced by $X_i$.

The expectation and covariance of the Gaussian process $Y_n^w$ are given by

$$\mu_n^w(t) = -(b_n^2 + w)(1 - C(s_n t, 0)), \tag{24}$$

$$r_n(t_1, t_2) = b_n^2(C(s_n t_1, s_n t_2) - C(s_n t_1, 0)C(s_n t_2, 0)). \tag{25}$$

Note that the conditional covariance $r_n(t_1, t_2)$ does not depend on $w$. Let $t, t_1, t_2 \in \mathbb{R}^d$, $w \in \mathbb{R}$ be fixed. Using (20) and (22), we obtain

$$\lim_{n \to \infty} \mu_n^w(t) = -\gamma(t), \tag{26}$$

$$\lim_{n \to \infty} r_n(t_1, t_2) = \gamma(t_1) + \gamma(t_2) - \gamma(t_1 - t_2). \tag{27}$$

A further consequence of (24) is that as long as $t$ stays bounded, there is a constant $c$ such that, for sufficiently large $n$, we have

$$|\mu_n^w(t)| \leq c + |w|/2 \qquad \forall w \in \mathbb{R}. \tag{28}$$

It follows from (26), (27) that as $n \to \infty$, the process $Y_n^w$ converges in the sense of finite-dimensional distributions to $\{W(t) - \gamma(t), t \in \mathbb{R}^d\}$, where $\{W(t), t \in \mathbb{R}^d\}$ is a Gaussian process with stationary increments, variogram $2\gamma$ and $W(0) = 0$. On the other hand, it is known (see, e.g., Corollary 4.19 in [26]), that the point process $\sum_{i=1}^n \delta_{Y_{i,n}(0)}$ converges, as $n \to \infty$, to the Poisson point process on $\mathbb{R}$ with intensity $e^{-y} dy$. From these two facts, at least on the formal level, we obtain the statement of the theorem. However, making this rigorous requires some work.

First, we show that $\eta_n$ converges to $\eta$ in the sense of finite-dimensional distributions. Let $t_1, \ldots, t_k \in \mathbb{R}^d$ and $y_1, \ldots, y_k \in \mathbb{R}$ be fixed. By conditioning on $Y_n(0) = w$ and noting that the density of $Y_n(0)$ is given by

$$f_{Y_n(0)}(w) = 1/(\sqrt{2\pi} b_n) e^{-(w + b_n^2)^2 / 2 b_n^2},$$

we obtain

$$\mathbb{P}[\exists j : Y_n(t_j) > y_j]$$
$$= \frac{1}{\sqrt{2\pi} b_n} \int_{\mathbb{R}} e^{-(w + b_n^2)^2 / (2 b_n^2)} \mathbb{P}[\exists j : Y_n(t_j) > y_j | Y_n(0) = w] \, dw$$
$$= \frac{1}{\sqrt{2\pi} b_n e^{b_n^2/2}} \int_{\mathbb{R}} e^{-w - w^2 / (2 b_n^2)} \mathbb{P}[\exists j : Y_n^w(t_j) > y_j - w] \, dw.$$

Noting that (21) implies that $\sqrt{2\pi} b_n e^{b_n^2/2} \sim n$ as $n \to \infty$ and taking $A > 0$, we may write the above as

$$\mathbb{P}[\exists j : Y_n(t_j) > y_j] \sim \frac{1}{n}\left(\int_{-A}^{A} + \int_{A}^{\infty} + \int_{-\infty}^{-A}\right) = \frac{1}{n}(I_1(n) + I_2(n) + I_3(n)).$$

Since the convergence of the distribution of $\{Y_n^w(y_j)\}_{j=1}^k$ to that of $\{W(y_j) - \gamma(y_j)\}_{j=1}^k$ is uniform provided that $w \in [-A, A]$, we obtain

$$(29) \qquad \lim_{n \to \infty} I_1(n) = \int_{-A}^{A} e^{-w} \mathbb{P}[\exists j : W(t_j) - \gamma(t_j) > y_j - w] \, dw.$$

For $I_2(n)$, we have the trivial estimate

$$(30) \qquad I_2(n) \leq \int_{A}^{\infty} e^{-w} \, dw = e^{-A}.$$

We estimate $I_3(n)$. Using (28), we obtain, if $w < -A$ and $A, n$ are large,

$$\mathbb{P}[Y_n^w(t_j) > y_j - w] \leq \mathbb{P}[Y_n^w(t_j) - \mu_n^w(t_j) > y_j - c - |w|/2 - w]$$
$$\leq \mathbb{P}[Y_n^w(t_j) - \mu_n^w(t_j) > |w|/4].$$

Recall the well-known estimate $\Psi(t) \leq e^{-t^2/2}$, $t \geq 0$, where $\Psi(t)$ is the tail of the standard Gaussian distribution. By (27), $\text{Var}[Y_n^w(t_j)] < \kappa^2$ for some $\kappa > 0$ and all $j = 1, \ldots, k$, $w \in \mathbb{R}$, $n \in \mathbb{N}$. Hence,

$$\mathbb{P}[Y_n^w(t_j) > y_j - w] \leq e^{-(w/4)^2/(2\kappa^2)}.$$

It follows that

$$I_3(n) \leq \sum_{j=1}^{k} \int_{-\infty}^{-A} e^{-w} \mathbb{P}[Y_n^w(t_j) > y_j - w] \, dw \leq k \int_{-\infty}^{-A} e^{-w} e^{-w^2/(32\kappa^2)} \, dw.$$

Hence,

$$(31) \qquad \lim_{A \to \infty} \limsup_{n \to \infty} I_3(n) = 0.$$

Bringing (29), (30) and (31) together and letting $A \to \infty$, we obtain

$$\mathbb{P}[\exists j : Y_n(t_j) > y_j] \sim \frac{1}{n} \int_{\mathbb{R}} e^{-w} \mathbb{P}[\exists j : W(t_j) - \gamma(t_j) > y_j - w] \, dw$$
$$= \frac{1}{n} \mathbb{E} \exp \max_{j=1,\ldots,k} (W(t_j) - \gamma(t_j) - y_j)$$

as $n \to \infty$. Therefore,

$$\lim_{n \to \infty} \mathbb{P}[\forall j : \eta_n(t_j) \leq y_j] = \lim_{n \to \infty} (1 - \mathbb{P}[\exists j : Y_n(t_j) > y_j])^n$$
$$= \exp\left\{-\mathbb{E} \exp \max_{j=1,\ldots,k} (W(t_j) - \gamma(t_j) - y_j)\right\}.$$

By (5), the right-hand side coincides with $\mathbb{P}[\forall j : \eta(t_j) \leq y_j]$, which proves that $\eta_n$ converges to $\eta$ in the sense of finite-dimensional distributions.

It remains to show that the sequence $\eta_n$ is tight in $C(K)$, where $K \subset \mathbb{R}^d$ is a fixed compact set. First, note that the sequence $\eta_n(0)$ is tight in $\mathbb{R}$ [in fact, the distribution of $\eta_n(0)$ converges weakly to the Gumbel distribution]. For a function $f \in C(K)$ and $\delta > 0$, define

$$\omega_\delta(f) = \sup_{t_1, t_2 \in K, \|t_1 - t_2\| \leq \delta} |f(t_1) - f(t_2)|.$$

By the standard tightness criterion (see, e.g., Theorem 7.3 in [3]), we need to show that for every $\varepsilon > 0$, $a > 0$, there exists some $\delta > 0$ such that

(32) $\qquad \mathbb{P}[\omega_\delta(\eta_n) > a] < \varepsilon \qquad$ for all $n > N$.

Throughout, $N$ denotes a large integer whose value may change from line to line. We concentrate on proving (32). The proof of the next lemma will be given later.

LEMMA 19. *The following assertions hold:*

1. *for every $c > 0$, the family of processes $Y_n^w$, $w \in [-c, c]$, $n \in \mathbb{N}$, is tight in $C(K)$;*
2. *the family of processes $Y_n^w - \mu_n^w$, $w \in \mathbb{R}$, $n \in \mathbb{N}$, is tight in $C(K)$.*

For $c_1 > 0$, define a sequence of random events

$$E_n = \left\{ \inf_{t \in K} \eta_n(t) < -c_1 \right\}.$$

We show that we can find $c_1 > 0$ such that $\mathbb{P}[E_n] < \varepsilon$ for all $n > N$. First, choose $c_0$ so large that $2e^{-c_0} < \varepsilon$. Using part 1 of Lemma 19, choose $c_1$ so large that

$$\mathbb{P}\left[ \inf_{t \in K} Y_n^w(t) < c_0 - c_1 \right] < 1/2 \qquad \text{for all } w \in [-c_0, c_0], n \in \mathbb{N}.$$

Define random events

$$A_{i,n} = \left\{ Y_{i,n}(0) \in [-c_0, c_0], \inf_{t \in K} Y_{i,n}(t) - Y_{i,n}(0) \geq c_0 - c_1 \right\}.$$

We have, by conditioning on $Y_{i,n}(0) = w$,

$$\mathbb{P}[A_{i,n}] = (\sqrt{2\pi} b_n e^{b_n^2/2})^{-1} \int_{-c_0}^{c_0} e^{-w - w^2/(2b_n^2)} \mathbb{P}\left[ \inf_{t \in K} Y_n^w(t) \geq c_0 - c_1 \right] dw$$

$$\geq \frac{1}{4n} \int_{-c_0}^{c_0} e^{-w - w^2/(2b_n^2)} dw, \qquad n > N,$$

which implies that $\mathbb{P}[A_{i,n}] \geq c_0/n$ if $c_0$ is sufficiently large and $n > N = N(c_0)$. Noting that $\mathbb{P}[E_n] \leq \mathbb{P}[\bigcap_{i=1}^{n} A_{i,n}^c]$ gives

$$\mathbb{P}[E_n] \leq (1 - c_0/n)^n \leq 2e^{-c_0} < \varepsilon, \qquad n > N.$$

For $c_2 > 0$, define the random events

$$F_n = \bigcup_{i=1}^{n} \{Y_{i,n}(0) > c_2\},$$

$$G_n = \Big\{\exists t \in K : \eta_n(t) \neq \sup_{i \in \{1,\ldots,n\} : |Y_{i,n}(0)| < c_2} Y_{i,n}(t)\Big\}.$$

Trivially, $\mathbb{P}[F_n] = \mathbb{P}[\eta_n(0) > c_2] < \varepsilon$ for every $n$, if $c_2$ is large. We show that there exists some $c_2 > 0$ such that $\mathbb{P}[G_n] < 3\varepsilon$ for $n > N$. Introduce random events

$$B_{i,n} = \Big\{Y_{i,n}(0) < -c_2, \sup_{t \in K} Y_{i,n}(t) > -c_1\Big\}.$$

Then, again conditioning on $Y_{i,n}(0) = w$ and recalling (28), we obtain

$$\mathbb{P}[B_{i,n}] = (\sqrt{2\pi} b_n e^{b_n^2/2})^{-1} \int_{-\infty}^{-c_2} e^{-w - w^2/(2b_n^2)} \mathbb{P}\Big[\sup_{t \in K} Y_n^w(t) > -c_1 - w\Big] dw.$$

By part 2 of Lemma 19, there exists some $c_3 > 0$ such that

$$\mathbb{P}\Big[\sup_{t \in K}(Y_n^w(t) - \mu_n^w(t)) > c_3\Big] < 1/2, \qquad w \in \mathbb{R}, n \in \mathbb{N}.$$

Recall that, by (28) and (25), we have

$$\sup_{t \in K} \mu_n^w \leq c_4 - \frac{w}{2}, \qquad \sup_{t \in K} \operatorname{Var} Y_n^w(t) < \kappa^2, \qquad w < 0, n > N,$$

for some $c_4, \kappa$. Applying Borell's inequality (see Theorem D.1 in [25]), together with the above estimates, we obtain, for $w < 0$,

$$\mathbb{P}\Big[\sup_{t \in K} Y_n^w(t) > -c_1 - w\Big] < 2\Psi(-(-c_1 - w/2 - c_3 - c_4)/\kappa),$$

where $\Psi$ is the tail of the standard Gaussian distribution. If $w < -4(c_1 + c_3 + c_4)$, this, together with the bound $\Psi(t) \leq e^{-t^2/2}$, $t \geq 0$, implies that

$$\mathbb{P}\Big[\sup_{t \in K} Y_n^w(t) > -c_1 - w\Big] < 2e^{-w^2/(32\kappa^2)}, \qquad n > N.$$

It follows that, for $n > N$ and $c_2 > 4(c_1 + c_3 + c_4)$,

$$\mathbb{P}[B_{i,n}] \leq \frac{4}{n} \int_{-\infty}^{-c_2} e^{-w} e^{-w^2/(32\kappa^2)} dw.$$

So, we can choose $c_2$ sufficiently large that $n\mathbb{P}[B_{1,n}] < \varepsilon$ for $n > N$. Therefore,

$$\mathbb{P}[G_n] \leq \mathbb{P}[E_n] + \mathbb{P}[F_n] + \mathbb{P}[G_n \setminus (E_n \cup F_n)] < 2\varepsilon + n\mathbb{P}[B_{1,n}] < 3\varepsilon.$$

We are now ready to prove (32). Let

$$C_{i,n} = \{Y_{i,n}(0) \in [-c_2, c_2], \omega_\delta(Y_{i,n}) > a\}$$

and define $H_n = \bigcup_{i=1}^n C_{i,n}$. Then,

$$\mathbb{P}[C_{i,n}] = (\sqrt{2\pi} b_n e^{b_n^2/2})^{-1} \int_{-c_2}^{c_2} e^{-w - w^2/(2b_n^2)} \mathbb{P}[\omega_\delta(Y_n^w) > a] \, dw.$$

By part 1 of Lemma 19 and the tightness criterion (see Theorem 7.3 in [3]), we can make $\mathbb{P}[\omega_\delta(Y_n^w) > a]$ arbitrary small (uniformly in $w \in [-c_2, c_2]$ and for $n > N$) by choosing $\delta$ small. So, choose $\delta > 0$ sufficiently small that $\mathbb{P}[C_{i,n}] < \frac{\varepsilon}{n}$. Then,

$$\mathbb{P}[\omega_\delta(\eta_n) > a] \leq \mathbb{P}[G_n] + \mathbb{P}[H_n] < 3\varepsilon + n\mathbb{P}[C_{1,n}] < 4\varepsilon,$$

which yields (32) with $4\varepsilon$ instead of $\varepsilon$. This proves the tightness of the sequence $\eta_n$ and completes the proof of Theorem 17. □

PROOF OF LEMMA 19. It follows from (25) that, independently of $w$,

$$\begin{aligned}\text{Var}(Y_n^w(t_1) - Y_n^w(t_2)) \\ = b_n^2 (2 - 2C(s_n t_1, s_n t_2) - (C(s_n t_1, 0) - C(s_n t_2, 0))^2) \\ \leq 2b_n^2 (1 - C(s_n t_1, s_n t_2)).\end{aligned}$$

Assumption 16 implies that, uniformly in $t_1, t_2 \in K$,

$$\text{Var}(Y_n^w(t_1) - Y_n^w(t_2)) \leq 2b_n^2 \cdot 2(L(s_n) s_n^\alpha \gamma(t_1 - t_2)), \qquad n > N.$$

By (22), we have $b_n^2 L(s_n) s_n^\alpha \leq 2$, $n > N$, and so, for some $c_5 > 0$,

$$(33) \quad \text{Var}(Y_n^w(t_1) - Y_n^w(t_2)) \leq 8\gamma(t_1 - t_2) \leq c_5 \|t_1 - t_2\|^\alpha, \qquad n > N.$$

Now, the second claim of the lemma follows from (33) by applying Corollary 11.7 of [22] to the family of processes $Y_n^w - \mu_n^w$, $w \in \mathbb{R}$, $n \in \mathbb{N}$ [take $\psi(x) = x^2$, $d^2(t_1, t_2) = c_5 \|t_1 - t_2\|^\alpha$ there]. To prove the first claim, we need to additionally show that $\mu_n^w$, $w \in [-c, c]$, $n \in \mathbb{N}$, is a tight family of functions in $C(K)$. This last statement follows from (24), which shows that the convergence $\mu_n^w(t) \to -\gamma(t)$ in (26) is uniform in $t \in K$, $w \in [-c, c]$. □

**7. Domains of attraction.** We are now going to prove a partial converse of Theorem 17. More precisely, we characterize all nontrivial limits of normalized and spatially rescaled pointwise maxima of stationary Gaussian processes. Let us call a random process $\{\eta(t), t \in \mathbb{R}^d\}$ *degenerate* if, for all $t_1, t_2 \in \mathbb{R}^d$, we have $\eta(t_1) = \eta(t_2)$ a.s.

THEOREM 20. *Let $\{X(t), t \in \mathbb{R}^d\}$ be a stationary zero-mean, unit-variance Gaussian process with continuous covariance $C(t) = \mathbb{E}[X(0)X(t)]$ and let $X_i$, $i \in \mathbb{N}$, be independent copies of $X$. Suppose that, for some sequences $a'_n > 0$, $b'_n \in \mathbb{R}$ and $s'_n > 0$, the process $\{\eta'_n(t), t \in \mathbb{R}^d\}$ defined by*

$$\eta'_n(t) = \bigvee_{i=1}^n a'_n(X_i(s'_n t) - b'_n)$$

*converges, as $n \to \infty$, to some nondegenerate, continuous-in-probability process $\{\eta'(t), t \in \mathbb{R}^d\}$, in the sense of finite-dimensional distributions. The following assertions then hold:*

1. *there is an $\alpha \in (0,2]$, a finite measure $\mu$ on the unit sphere $\mathbb{S}^{d-1}$ in $\mathbb{R}^d$ and a function $L$ that varies slowly at 0 such that*

(34) $$1 - C(t) \sim L(\|t\|)\gamma(t) \quad \text{as } t \to 0,$$

   *where*

(35) $$\gamma(t) = \int_{\mathbb{S}^{d-1}} |\langle t, x \rangle|^\alpha \, d\mu(x);$$

2. *the normalizing sequences $a'_n, b'_n$ and $s'_n$ satisfy*

(36) $$\lim_{n \to \infty} a'_n/b_n = A > 0, \quad \lim_{n \to \infty} b_n(b'_n - b_n) = B \in \mathbb{R},$$

(37) $$\lim_{n \to \infty} b_n^2 L(s'_n) s'^\alpha_n = s > 0,$$

   *where $b_n$ is defined by (21);*
3. *the limiting process $\eta'$ coincides with $A(\eta - B)$, where $\eta$ is the Brown–Resnick process associated to the variogram $2s\gamma$.*

We need a lemma, the essential part of which was proven in [18].

LEMMA 21. *For $n \in \mathbb{N}$, let $Z_1^{(n)}, \ldots, Z_n^{(n)}$ be i.i.d. bivariate Gaussian vectors having standard Gaussian margins and correlation $\rho_n$. The maxima*

$$M_n = \bigvee_{i=1}^n b_n(Z_i^{(n)} - b_n)$$

*converge in distribution to some bivariate random vector if and only if*

(38) $$\lim_{n \to \infty} b_n^2(1 - \rho_n) = c$$

*for some* $c \in [0, \infty]$. *The limiting bivariate distribution depends on* $c$ *continuously; its margins are independent if and only if* $c = \infty$ *and are equal a.s. if and only if* $c = 0$.

PROOF. Suppose, first, that (38) holds. Then, by a result of [18], the sequence $M_n$ converges in distribution. The explicit formula, given in [18], shows that the limiting distributions corresponding to different values of $c$ are different. Suppose, now, that (38) does not hold. We then have $0 \leq \liminf b_n^2(1 - \rho_n) < \limsup b_n^2(1 - \rho_n) \leq \infty$. Again using [18], we obtain that the sequence $M_n$ has at least two different accumulation points and thus does not converge. The last claim of the lemma follows again from the explicit formula in [18]. $\square$

PROOF OF THEOREM 20. By stationarity of $X$, the distribution of $\eta'(t)$ does not depend on $t \in \mathbb{R}^d$. Thus, if for some constant $c_0$, $\eta'(0) = c_0$ a.s., then for every $t \in \mathbb{R}^d$, $\eta'(t) = c_0$ a.s., which is a contradiction since $\eta'$ is assumed to be nondegenerate. So, in the sequel, we assume that $\eta'(0)$ is not a.s. constant. In this case, the convergence-to-types theorem (see Proposition 0.2 in [26]), together with (23), yields constants $A > 0$, $B \in \mathbb{R}$ such that (36) holds. It follows that the process

$$\eta_n(t) = \bigvee_{i=1}^{n} b_n(X_i(s'_n t) - b_n)$$

converges, as $n \to \infty$, to the nondegenerate limit $\eta = A^{-1}\eta' + B$. From now on, we consider the processes $\eta_n$ and $\eta$ instead of $\eta'_n$ and $\eta'$.

For any fixed $t \in \mathbb{R}^d$, the previous lemma, applied to the triangular array of bivariate vectors $Z_i^{(n)} = (X_i(0), X_i(s'_n t))$, $i = 1, \ldots, n$, $n \in \mathbb{N}$, yields a constant $c(t) \in [0, \infty]$ such that

(39) $$\lim_{n \to \infty} b_n^2(1 - C(s'_n t)) = c(t).$$

Since the limiting process $\eta$ is assumed to be continuous in probability, the distribution of the bivariate vector $(\eta(0), \eta(t))$ must converge weakly to the distribution of $(\eta(0), \eta(0))$ as $t \to 0$. Using the last statement of Lemma 21, we obtain that $\lim_{t \to 0} c(t) = c(0) = 0$, that is, $c$ is continuous at the origin. Note, also, that $c(t_0) \neq 0$ for some $t_0 \neq 0$ since otherwise the process $\eta$ would be degenerate.

By Bochner's theorem, there exists an $\mathbb{R}^d$-valued random variable $\xi$ such that the characteristic function of $\xi$ is $C(t)$. Moreover, since the function $C$ is real-valued, the distribution of $\xi$ must be symmetric with respect to the origin. Let $\xi_i$, $i \in \mathbb{N}$, be i.i.d. copies of $\xi$. Then, the characteristic function $\varphi_n$ of

$$S_n = s'_n \sum_{i=1}^{[b_n^2]} \xi_i$$

is given by $\varphi_n(t) = C(s'_n t)^{[b_n^2]}$ so that (39) yields

$$\lim_{n\to\infty} \varphi_n(t) = \lim_{n\to\infty} (1 - c(t)/b_n^2 + o(1/b_n^2))^{[b_n^2]} = e^{-c(t)}.$$

Now, Lévy's convergence theorem tells us that the random vector $S_n$ converges weakly to a random vector $S$ whose distribution is necessarily nondegenerate (i.e., $\mathbb{P}[S=0] \neq 1$; to see this, recall that $c(t_0) \neq 0$ and hence $e^{-c(t_0)} \neq 1$ for some $t_0 \neq 0$), stable with some parameter $\alpha \in (0,2]$ and symmetric with respect to the origin. It follows from the characterization of domains of attraction of multidimensional symmetric stable distributions in terms of characteristic functions (see Corollaries 1 and 2 in [1]) that the covariance function $C$ must have the form (34), (35). Inserting this in (39) for some $t$ with $\|t\| = 1$, we obtain

$$\lim_{n\to\infty} b_n^2 L(s'_n) s'^\alpha_n \gamma(t) = c(t),$$

which yields (37). Furthermore, (34) and (35) imply that the process $X$ satisfies Assumption 16. Therefore, by Theorem 17, the limiting process $\eta$ must be the Brown–Resnick process associated to the variogram $2s\gamma$. Recalling that $\eta = A^{-1}\eta' + B$, we obtain the last statement of the theorem. $\square$

**8. Extensions and remarks.** In view of Theorems 17 and 20, the question arises as to whether max-stable processes corresponding to variograms $\gamma$ that are not of the form (35) also admit representations as limits of pointwise maxima of stationary Gaussian processes in some broader sense, as in Theorem 20. The answer is affirmative, as the following theorem shows.

THEOREM 22. *Let $\gamma$ be a variogram on $\mathbb{R}^d$, that is, $\gamma(0) = 0$ and $\gamma$ is negative definite. For each $n \in \mathbb{N}$, let $X_{1n}, \ldots, X_{nn}$ be i.i.d. copies of a stationary zero-mean Gaussian process $\{X_n(t), t \in \mathbb{R}^d\}$ with covariance function $\exp(-\gamma(t)/b_n^2)$. Define*

$$\eta_n(t) = \bigvee_{i=1}^n b_n(X_{in}(t) - b_n), \qquad t \in \mathbb{R}^d.$$

*Then, $\eta_n$ converges, in the sense of finite-dimensional distributions, to the Brown–Resnick process associated to the variogram $2\gamma$.*

PROOF. Note that $\exp(-\gamma(t)/b_n^2)$ is indeed a covariance function of some stationary Gaussian process, by Schoenberg's theorem (see Theorem 7.8 in [2]). As in the proof of Theorem 17, it can be shown that the conditional distribution of $b_n(X_{in}(t) - X_{in}(0))$, given that $b_n(X_{in}(0) - b_n) = w$, converges to the distribution of $W(t) - \gamma(t)$, where $W$ is a Gaussian process with stationary increments, variogram $2\gamma$ and $W(0) = 0$. The rest of the proof is the same as that of Theorem 17. $\square$

REMARK 23. The above theorem gives another proof of stationarity in Theorem 2 in the case $W(0) = 0$.

REMARK 24. The bivariate distributions of the Brown–Resnick process $\eta$ associated to the variogram $\gamma$ are given by the formula

$$\mathbb{P}(\eta(t_1) \leq y_1, \eta(t_2) \leq y_2)$$
$$= \exp\left\{-e^{-y_1}\Phi\left(\sqrt{\gamma(t_1-t_2)}/2 + \frac{y_2-y_1}{\sqrt{\gamma(t_1-t_2)}}\right)\right.$$
$$\left. - e^{-y_2}\Phi\left(\sqrt{\gamma(t_1-t_2)}/2 + \frac{y_1-y_2}{\sqrt{\gamma(t_1-t_2)}}\right)\right\},$$

where $\Phi$ is the standard normal distribution function.

PROOF. The remark is a consequence of Theorem 22 and a result of [18]. Moreover, it follows from Theorem 22 that the finite-dimensional distributions of the process $\eta$ belong to the family of multivariate max-stable distributions introduced in [18]. □

REMARK 25. A natural dependence measure between $\eta(0)$ and $\eta(t)$ is given by $\rho(t) = 2 - \varsigma(t) \in [0,1]$, where $\varsigma(t)$ is determined from the condition

$$\mathbb{P}[\eta(0) \leq z, \eta(t) \leq z] = \mathbb{P}[\eta(0) \leq z]^{\varsigma(t)}$$

for some (and hence all) $z \in \mathbb{R}$; see, for example, [7, 28]. It follows from Remark 24 that

$$\rho(t) = 2(1 - \Phi(\sqrt{\gamma(t)}/2)).$$

Thus, a variogram $\gamma$ is completely determined by the dependence function $\rho(t)$ of the corresponding process $\eta$. It follows that $\eta(0)$ and $\eta(t)$ become asymptotically independent as $\|t\| \to \infty$ [which corresponds to $\rho(t) \to 0$] if and only if $\gamma(t) \to \infty$ as $\|t\| \to \infty$. Furthermore, if $d = 1$, then, by Theorem 3.4 in [30], the process $\eta$ is mixing if and only if $\gamma(t) \to \infty$ as $t \to \infty$.

REMARK 26. Theorem 17 may be generalized to processes whose covariance has different Hölder exponents in different directions. For example, assume that $\{X(t), t \in \mathbb{R}^d\}$ is a stationary zero-mean Gaussian process with covariance function $C$ satisfying

$$C(t) = C(t_1, \ldots, t_d) = 1 - \sum_{i=1}^{d} c_i |t_i|^{\alpha_i} + o(\|t\|^{\alpha_d}) \qquad \text{as } t \to 0$$

for some $0 < \alpha_1 \leq \cdots \leq \alpha_d \leq 2$, $c_1, \ldots, c_d > 0$. If $X_i, i \in \mathbb{N}$, are independent copies of $X$, then

$$\eta_n(t) = \bigvee_{i=1}^n b_n(X_i(b_n^{-2/\alpha_1} t_1, \ldots, b_n^{-2/\alpha_d} t_d) - b_n)$$

converges to the Brown–Resnick process associated to the variogram $2\gamma$, where $\gamma(t_1, \ldots, t_d) = \sum_{i=1}^d c_i |t_i|^{\alpha_i}$.

**Acknowledgments.** The authors are grateful to the referees for numerous useful suggestions which considerably improved the paper. We thank also Achim Wübker for introducing us to [1].

Z. KABLUCHKO
M. SCHLATHER
INSTITUT FÜR MATHEMATISCHE STOCHASTIK
GEORG-AUGUST-UNIVERSITÄT GÖTTINGEN
GOLDSCHMIDTSTR. 7 D-37077 GÖTTINGEN
GERMANY
E-MAIL: kabluch@math.uni-goettingen.de
       schlather@math.uni-goettingen.de

L. DE HAAN
DEPARTMENT OF ECONOMICS
ERASMUS UNIVERSITY ROTTERDAM
P.O. BOX 1738
3000 DR, ROTTERDAM
THE NETHERLANDS
E-MAIL: ldhaan@few.eur.nl